\input eplain

\baselineskip=14pt
\parskip=10pt

\magnification=\magstephalf

\def\1{{\overline{1}}}
\def\2{{\overline{2}}}
\parindent=0pt
\overfullrule=0in

\def\frac#1#2{{#1 \over #2}}

\bf
\centerline
{
The Number of 1...d-Avoiding Permutations of Length d+r for SYMBOLIC d but Numeric r 
}
\rm
\bigskip
\centerline
{\it By Shalosh B. EKHAD, Nathaniel SHAR, and Doron ZEILBERGER}
\bigskip
\qquad \qquad 
{\it Dedicated to Ira Martin GESSEL (b. April 9, 1951), on his millionth$_2$ birthday}
\bigskip

{\bf Preface: How many permutations are there of length googol+30 avoiding an increasing subsequence of length googol?}

This number is way too big for our physical universe, but the number of permutations 
of length googol+30 that {\it contain} at least one increasing subsequence of length googol is
a certain integer that may be viewed in
{\tt http://www.math.rutgers.edu/\~{}zeilberg/mamarim/mamarimPDF/gessel64.pdf}.
Hence the number of permutations of length googol+$30$ {\it avoiding} an increasing subsequence of length googol is
$(googol+30)!$ {\it minus} the above small number.

{\bf Counting the ``Bad Guys''}

Recall that thanks to Robinson-Schensted ([Rob][Sc]), the number of permutations of length $n$ that 
do {\bf not} contain an increasing subsequence of length $d$ is given by
$$
G_d(n):=\sum_{ {\lambda \vdash n} \atop {\#rows(\lambda)<d} } f_{\lambda}^2 \quad ,
$$
where $\lambda$ denotes a typical {\it Young diagram}, and $f_\lambda$ is the number of
{\it Standard Young tableaux} whose {\it shape} is $\lambda$.

Hence the number of permutations of length $n$ that {\bf do} contain an increasing subsequence of length $d$ is
$$
B_d(n):=\sum_{ {\lambda \vdash n} \atop {\#rows(\lambda) \geq d} } f_{\lambda}^2 \quad .
$$
Since the total number of permutations  of length $n$ is $n!$ ([B]), if we know how to find $B_d(n)$, we would
know immediately $G_d(n)=n!-B_d(n)$, at least if we leave $n!$ alone as a factorial, rather than spell it out.

Recall that the {\it Hook Length formula} (see [Wiki]) tells you that if $\lambda$ is a Young diagram
then
$$
f_{\lambda}= \frac{n!}{\prod_{c \in \lambda} h(c)} \quad,
$$
where the product is over all the $n$ cells of the Young diagram, and the {\it hook-lenght}, $h(c)$, of
a cell $c=(i,j)$, is $(\lambda_i -i)+ (\lambda'_j -j)+1$, where $\lambda'$ is the {\it conjugate} diagram, where
the rows become columns and vice-versa.

Let $r$ be a fixed integer, then for {\it symbolic} $d$, valid for $d \geq r-1$, any
Young diagram with at least $d$ rows, and with $d+r$ cells, can be written, for some Young diagram 
$\mu=(\mu_1, \dots, \mu_r)$, with $\leq r$ cells,
(where we add zeros to the end if the number of parts of $\mu$ is less than $r$)
as
$$
\lambda=(1+\mu_1, \dots, 1+\mu_r, 1^{d-r+r'} )\quad ,
$$
where $r'=r-|\mu|$. For such a shape $\lambda$, with {\it at least} $d$ rows,
$$
\prod_{c \in \lambda} h(c)= \left ( \, \prod_{c \in \mu} h(c) \, \right ) \cdot ( (d+r'+\mu_1)(d+r'-1+\mu_2) \cdots (d+r'-r+1+\mu_r) ) \cdot (d-r+r')!  \quad .
$$
Hence $f_{\lambda}$, that is $(d+r)!$ divided by the above, is a certain specific number times a certain polynomial in $d$.
Since, for a specific, {\it numeric},  $r$, there are only {\it finitely} many Young diagrams with at most $r$ cells, the
computer can find all of them, compute the polynomial corresponding to each of them, square it, and add-up all these terms, getting
an {\it explicit} {\bf polynomial} 
expression, in the variable $d$, for $B_d(d+r)$, the number of permutations of length $d+r$ that {\it contain}
an increasing subsequence of length $d$. As we said above, from this we can find $G_d(d+r)=(d+r)!-B_d(d+r)$, 
valid for {\it symbolic} $d \geq r-1$.

{\bf ${\bf B_d(d+r)}$ for r from 0 to 30}
$$
B_d(d)= 1 \quad ,
$$
$$
B_d(d+1)= {d}^{2}+1 \quad ,
$$
$$
B_d(d+2)= \frac{1}{2} \, {d}^{4}+{d}^{3}+  \frac{1}{2}\,{d}^{2}+d+3 \quad ,
$$
$$
B_d(d+3)=  \frac{1}{6} \,{d}^{6}+{d}^{5}+  \frac{5}{3}\,{d}^{4}+\frac{2}{3}\,{d}^{3}+{\frac {19}{6}}\,{d}^{2}+{\frac {31}{3}}\,d+11 \quad ,
$$
$$
B_d(d+4)= 
\frac{1}{24}\,{d}^{8}
+\frac{1}{2}\,{d}^{7}
+{\frac {25}{12}}\,{d}^{6}
+{\frac {19}{6}}\,{d}^{5} 
+{\frac {29}{24}}\,{d}^{4}
+9\,{d}^{3}
+{\frac {247}{6}}\,{d}^{2}
+{\frac {395}{6}}\,d
+47
\quad ,
$$
$$
B_d(d+5)= 
{\frac {1}{120}}\,{d}^{10}
+  \frac{1}{6}\,{d}^{9}
+{\frac {31}{24}}\,{d}^{8}
+  \frac{14}{3}\,{d}^{7}
+{\frac {823}{120}}\,{d}^{6}
+{\frac {67}{30}}\,{d}^{5} 
+{\frac {653}{24}}\,{d}^{4}
+{\frac {959}{6}}\,{d}^{3}
+{\frac {10459}{30}}\,{d}^{2}
+{\frac {3981}{10}}\,d
+239 \quad .
$$
For $B_d(d+r)$ for $r$ from $6$ up to $30$, see \hfill\break 
{\tt http://www.math.rutgers.edu/\~{}zeilberg/tokhniot/oGessel64a} .

{\bf Sequences}

The sequence $G_3(n)$ is the greatest {\it celeb} in the kingdom of combinatorial sequences
[the subject of an entire book([St]) by Ira Gessel's 
illustrious {\it academic father}, Richard Stanley],
the super-famous {\bf A000108} in 
Neil Sloane's legendary database ([Sl]). $G_4(n)$, while not in the same league as 
the Catalan sequence,  is still moderately famous,
{\bf A005802}. $G_5(n)$ is {\bf A047889}, 
$G_6(n)$ is {\bf A047890},
$G_7(n)$ is {\bf A052399 },
$G_8(n)$ is {\bf  A072131},
$G_9(n)$ is {\bf  A072132},
$G_{10}(n)$ is {\bf  A072133},
$G_{11}(n)$ is {\bf  A072167},
but $G_d(n)$ for $d \geq 12$ are absent (for a good reason, one must stop somewhere!).
Also the {\it flattened version} of the 
{\it double-sequence}, $\{G_{d}(n)\}$, for $1 \leq d \leq n \leq 45$ is {\bf A047887}.
Using the polynomials $B_d(d+r)$, 
we  computed the first $2d+1$ terms of $G_d(n)$ for $d \leq 30$. See \hfill\break 
{\tt http://www.math.rutgers.edu/\~{}zeilberg/tokhniot/oGessel64b} .

But this method can only go up to $2d+1$ terms of the sequence $G_d(n)$, and of course, the first
$d-1$ terms are trivial, namely $n!$. Can we find the first $100$ terms (or whatever) for the sequences
$G_d(n)$ for $d$ up to $20$, and beyond, {\bf efficiently?}

{\bf Encore: Efficient Computer-Algebra Implementation of Ira Gessel's AMAZING  Determinant Formula}

Recall Ira Gessel's [G]
famous expression for the generating function of $G_d(n)/n!^2$, {\it canonized} in the {\it bible} ([W], p. 996, Eq. (5)).
Here it is: 
$$
\sum_{n \geq 0} \frac{G_d(n)}{n!^2} \, x^{2n} = \, \det (I_{|i-j|}(2x))_{i,j=1, \dots, d}   \quad ,
$$
in which $I_{\nu}(t)$ is (the modified Bessel function)
$$
I_\nu (t) = \sum_{j=0}^{\infty} \frac{ (\frac{1}{2} \, t)^{2j+\nu}} {j!(j+\nu)!} \quad .
$$
Can we use this to compute the first $100$ terms of, say, $G_{20}(n)$?

While computing {\it numerical} determinants is very fast, computing {\it symbolic} ones is a different story.
First, do not get scared by the ``infinite'' power series. If we are only interested in the first $N$ terms of
$G_d(n)$, then it is safe to truncate the series up to $t^{2N}$, and take the determinant of a $d \times d$
matrix with {\it polynomial entries}. If you use the {\it vanilla} determinant in a computer-algebra
system such as Maple, it would be very inefficient, since the degree of the determinant is much  larger
than $2N$. But a little cleverness can make things more efficient. The Maple package
{\tt Gessel64}, available free of charge from

{\tt http://www.math.rutgers.edu/\~{}zeilberg/tokhniot/Gessel64} \quad ,

accompanying this article, has a procedure {\tt SeqIra(k,N)} that
computes the first N terms of $G_k(n)$, using a division-free
algorithm (see [Rot]) over an appropriate ring to compute the
determinant in Gessel's famous formula.

\verbatim
SeqIra:=proc(k,N) local ira,t,i,j, R:
  R := table():
  R[`0`] := 0:
  R[`1`] := 1:
  R[`+`] := `+`:
  R[`-`] := `-`:
  R[`*`] := proc(p, q): return add(coeff(p*q, t, i)*t**i, i=0..2*N): end:
  R[`=`] := proc(p, q): return evalb(p = q): end:
  ira:=expand(LinearAlgebra[Generic][Determinant][R](Matrix([seq([seq(Iv(abs(i-j),t,2*N),
                                                                        j=1..k-1)],
                                                                   i=1..k-1)]))):
  [seq(coeff(ira,t,2*i)*i!**2,i=1..N)]:
end:
|endverbatim

In the above code, procedure {\tt Iv(v,t,N)} computes the truncated modified Bessel function that shows
up in Gessel's determinant, and it is short enough to reproduce here:

{\tt Iv:=proc(v,t,N) local j: add(t**(2*j+v)/j!/(j+v)!,j=0..trunc((N-v)/2)+1): end: } \quad .

Using this procedure, the first-named author computed (in $4507$ seconds) the first $100$ terms of each of the sequences
$G_d(n)$ for $3 \leq d \leq 20$, and could have gone much further.

See {\tt http://www.math.rutgers.edu/\~{}zeilberg/tokhniot/oGessel64c} \quad .

{\bf HAPPY 64th BIRTHDAY, IRA!}

{\bf References}

[B] Rabbi Levi Ben Gerson, {\it Sefer Maaseh Hoshev}, Avignon, 1321.

[G] I.  Gessel, {\it Symmetric functions and P-recursiveness},
Journal of Combinatorial Theory Series A {\bf 53} (1990), 257-285; \quad 
{\tt http://people.brandeis.edu/\~{}gessel/homepage/papers/dfin.pdf} \quad .

[Rob] G. de B. Robinson, {\it On the representations of $S_n$}, Amer. J. Math. {\bf 60} (1938), 745-760.

[Rot] G. Rote, {\it Division-Free Algorithms for the Determinant and Pfaffian: Algebraic and Combinatorial Approaches}, 
{\it Computational Discrete Mathematics: Advanced Lectures} (2001), 119-135.

[Sc] C. E. Schensted, {\it Largest increasing and decreasing subsequences}, Canad. J. Math {\bf 13} (1961), 179-191.

[Sl] N.J.A. Sloane, {\it The On-Line Encyclopedia of Integer Sequences}; {\tt http://oeis.org} \quad .

[St] R. P. Stanley, {\it ``Catalan Numbers''}, Cambridge University Press, 2015.

[Wiki] The Wikipedia Foundation, {\it Hook Length Formula}; \hfill\break
{\tt http://en.wikipedia.org/wiki/Hook\_length\_formula}  \quad .

[W] H. S. Wilf,  {\it Mathematics, an experimental science}, in: 
``Princeton Companion to Mathematics'', (W. Timothy Gowers, ed.), Princeton University Press, 2008, 991-1000; \hfill\break
{\tt http://www.math.rutgers.edu/\~{}zeilberg/akherim/HerbMasterpieceEM.pdf} \quad .

\bigskip
\hrule
\bigskip
Shalosh B. Ekhad, c/o D. Zeilberger, Department of Mathematics, Rutgers University (New Brunswick), Hill Center-Busch Campus, 110 Frelinghuysen
Rd., Piscataway, NJ 08854-8019, USA.
\bigskip
\hrule
\bigskip
Nathaniel Shar, Department of Mathematics, Rutgers University (New Brunswick), Hill Center-Busch Campus, 110 Frelinghuysen
Rd., Piscataway, NJ 08854-8019, USA. \hfill \break
nshar at math dot rutgers dot edu \quad ;  \quad {\tt http://www.math.rutgers.edu/\~{}nbs48/} \quad .
\bigskip
\hrule
\bigskip
Doron Zeilberger, Department of Mathematics, Rutgers University (New Brunswick), Hill Center-Busch Campus, 110 Frelinghuysen
Rd., Piscataway, NJ 08854-8019, USA. \hfill \break
zeilberg at math dot rutgers dot edu \quad ;  \quad {\tt http://www.math.rutgers.edu/\~{}zeilberg/} \quad .

\bigskip
\hrule
\bigskip
Published in The Personal Journal of Shalosh B. Ekhad and Doron Zeilberger  \hfill \break
({ \tt http://www.math.rutgers.edu/\~{}zeilberg/pj.html})
and {\tt arxiv.org}. 
\bigskip
\hrule
\bigskip
{\bf April 9, 2015}

\end